\def\a{{\alpha}}

\def\c{{\gamma}}
\def\d{{\delta}}
\def\e{{\epsilon}}
\def\k{{\kappa}}
\def\t{{\tau}}
\def\T{{\sf t}}
\def\Th{{\Theta}}
\def\w{{\omega}}
\def\W{{\Omega}}

\def\bw{{\bf w}}

\def\bG{{\mathbf G}}

\def\sG{{\sf G}}
\def\sH{{\sf H}}

\def\N{{\mathbb N}}
\def\Z{{\mathbb Z}}

\def\B{{\cal B}}

\def\F{{\cal F}}
\def\G{{\cal G}}
\def\cH{{\cal H}}
\def\K{{\cal K}}
\def\P{{\cal P}}
\def\Q{{\cal Q}}

\newcommand{\bmp}{\mathcal{B}_{m,p}}
\newcommand{\gnp}{\mathcal{G}_{n,p}}
\newcommand{\cgQ}{\mathcal{P}}

\def\diam{{\sf diam}}
\def\dist{{\sf dist}}
\def\dom{{\sf dom}}
\def\gir{{\sf girth}}
\newcommand{\lgn}{\lg n}
\def\mt{{\emptyset}}
\def\lgN{{\lg N}}
\def\sub{{\ \subseteq\ }}

\def\pr{{\prime}}
\def\rar{{\rightarrow}}
\def\slN{{\sqrt\lgN}}

\def\sqr#1#2{{\vcenter{\hrule height.#2pt
        \hbox{\vrule width.#2pt height#1pt \kern#1pt
                \vrule width.#2pt}
        \hrule height.#2pt}}}
\def\gbox{{\mathchoice\sqr34\sqr34\sqr{2.1}3\sqr{1.5}3}}
\def\la{{\langle}}
\def\ra{{\rangle}}

\documentclass[12pt]{article}
\newtheorem{theorem}{Theorem}[section]

\newtheorem{result}[theorem]{Theorem}
\newtheorem{conjecture}[theorem]{Conjecture}

\newtheorem{question}[theorem]{Question}

\newtheorem{problem}[theorem]{Problem}

\usepackage{latexsym,amsmath,amsfonts,amssymb,graphicx}

\begin{document}

%##############################################################################
%##############################################################################
%
%       TITLE PAGE:
%
\title{Recent Progress in Graph Pebbling}

\author{
Glenn Hurlbert\thanks{Partially supported by National Security Agency
        grant \#H98230-05-1-0253.}~\\
Department of Mathematics and Statistics\\
Arizona State University\\
Tempe, Arizona 85287-1804\\
email: hurlbert@asu.edu\\
}
\maketitle
\newpage

%#############################################################################
%#############################################################################
%
%       ABSTRACT:
%
\begin{abstract}
The subject of graph pebbling has seen dramatic growth recently,
both in the number of publications and in the breadth of variations
and applications.
Here we update the reader on the many developments that have occurred
since the original {\it Survey of Graph Pebbling} in 1999.
\end{abstract}

\newpage

% ##########################################################################
%
% 		Introduction
%

\section{Introduction} \label{sec:intro}

Since the publication of the graph pebbling survey \cite{H} there has been
a great deal of activity in the subject.
There are now over 50 papers by roughly 80 authors in the field; 
the web page \cite{GPP} maintains a current list of these papers.
Many researchers have asked for an updated survey --- thanks to the 
New York Academy of Sciences for providing the opportunity.
The following is based on the talk, ``{\it Everything you always
wanted to know about graph pebbling but were afraid to ask}'', an obvious
ripoff/homage to one of New York's favorite directors.\footnote{Google 
``everything you always wanted to know'' if you do not understand this 
reference.}

We begin by introducing relevant terminology and background on the subject.
Here, the term {\it graph} refers to a simple graph without loops or 
multiple edges.
For the definitions of other graph theoretical terms see any standard
graph theory text such as \cite{W}.
The pebbling number of a disconnected graph will be seen to be undefined.
Henceforth we will assume all graphs to be connected.

A configuration $C$ of pebbles on a graph $G=(V,E)$ can be thought of 
as a function $C:V\rar\N$.
The value $C(v)$ equals the number of pebbles placed at vertex $v$,
and the {\it size} of the configuration is the number $|C|=\sum_{v\in V}C(v)$
of pebbles placed in total on $G$.
A pebbling step along an edge from $u$ to $v$ reduces by 2 the number 
of pebbles at $u$ and increases by 1 the number of pebbles at $v$.  
We say that a vertex $w$ can be {\it reached} by $C$ if one can repeatedly 
apply pebbling steps so that, in the resulting configuration $C^\pr$,
we have $C^\pr(w)\ge 1$ (and $C^\pr(v)\ge 0$ for all $v$).
The {\it pebbling number}, $\pi(G)$, is defined to be the smallest integer $m$ 
so that any specified {\it root} vertex $r$ of $G$ can be reached by every 
configuration $C$ of size $m$.
A configuration that reaches every vertex is called {\it solvable}, and
{\it unsolvable} otherwise.

% ##########################################################################
%
% 		Group Theory
%

\section{Group Theoretic Origins}\label{GroupTh}

The origins of graph pebbling reside in combinatorial number theory and
group theory.
A sequence of elements of a finite group $\sG$ is called a {\it zero-sum 
sequence} if it sums to the identity of $\sG$.
A simple pigeonhole argument (on the sequence of partial sums) proves the
following theorem.

\begin{theorem}\label{origin}
Any sequence of $|\sG|$ elements of a finite group $\sG$ contains a 
zero-sum subsequence.
\end{theorem}
In fact, a subsequence of consecutive terms can be guaranteed by the 
pigeonhole argument.
Furthermore, one can instead stipulate that the zero-sum subsequence has 
at most $N$ terms, where $N=N(\sG)$ is the exponent of $\sG$ (i.e. the maximum
order of an element of $\sG$), and this is best possible.

Initiated in 1956 by Erd\H{o}s \cite{Erd}, the study of zero-sum sequences 
has a long history with many important applications in number theory and 
group theory.
In 1961 Erd\H{o}s et al. \cite{EGZ} proved that every sequence of $2|\sG|-1$
elements of a cyclic group $\sG$ contains a zero-sum subsequence of length
exactly $|\sG|$.
In 1969 van Emde Boas and Kruyswijk \cite{EBK} proved that any sequence of
$N(1+\log(|\sG|/N))$ elements of a finite abelian group contains a zero-sum
sequence.
In 1994 Alford et al. \cite{AGP} used this result and modified Erd\H{o}s's
arguments to prove that there are infinitely many Carmichael numbers.
Much of the recent study has involved finding Davenport's constant $D(\sG)$,
defined to be the smallest $D$ such that every sequence of $D$
elements contains a zero-sum subsequence \cite{O}.
There are a wealth of results on this problem
\cite{Car,Gao,GG,GT,GS,Sun} and its variations \cite{GJ,Nat},
as well as applications to factorization theory \cite{Cha} and to graph 
theory \cite{AFK}.

In 1989 Kleitman and Lemke \cite{KL}, and independently Chung \cite{Chung},
proved the following theorem (originally stated number-theoretically),
strengthening Theorem \ref{origin}.
Let $\Z_n$ denote the cyclic group of order $n$, and let $|g|$ denote
the order of an element $g$ in the group to which it belongs.

\begin{result}\label{cyclic}
For every sequence $(g_k)_{k=1}^n$ of $n$ elements from $\Z_n$ there is 
a zero-sum subsequence $(g_k)_{k\in K}$ such that $\sum_{k\in K}1/|g_k|\le 1$.
\end{result}
For a sequence $S$ the sum $\sum_{g\in S}1/|g|$ is known as the {\it cross
number} of $S$ and is an important invariant in factorization theory.
Guaranteeing cross number at most 1 strengthens the extension of Theorem
\ref{origin} that $|K|\le N(\sG)$, and shows that equality holds if and only
if every $|g_k|=N$.  
The concept of pebbling in graphs arose from an attempt by Lagarias and
Saks to give an alternative (and more natural and structural) proof than 
that of Kleitman and Lemke; it was Chung who carried out their idea.
See also \cite{Denley} for another extension of this result.

Kleitman and Lemke then conjectured that Theorem \ref{cyclic} holds for all
finite groups.
For a subgroup $\sH$ of $\sG$, call a sequence of elements of $\sG$ an
$\sH$-{\it sum sequence} if its elements sum to an element of $\sH$.
In \cite{EH,G} is proved the following theorem (the methods of \cite{EH}
use graph pebbling).

\begin{theorem}\label{SubgroupConj}
Let $\sH$ be a subgroup of a finite abelian group $\sG$ with $|\sG|/|\sH|=n$.
For every sequence $(g_k)_{k=1}^n$ of $n$ elements from $\sG$ there is 
an $\sH$-sum subsequence $(g_k)_{k\in K}$ such that
$\sum_{k\in K}1/|g_k|\le 1/|\sum_{k\in K}\ g_k|$.
\end{theorem}
The case $\sH=\{e\}$ here gives Theorem \ref{cyclic} for finite abelian groups,
strengthening the van Emde Boas and Kruyswijk result \cite{EBK}.
Kleitman and Lemke also conjectured that Theorem \ref{SubgroupConj} holds for
all finite groups, and verified their conjecture for all dihedral groups 
(see \cite{KL}).
For other nonabelian groups, it has been shown recently to hold for the 
nonabelian solvable group of order 21 (see \cite{EH}).

It would be interesting to see whether graph pebbling methods can shed light
on the Davenport constant for finite abelian groups.
In this regard, one of the most pressing questions is as follows.
Write $\sG=\prod_{i=1}^r\Z_{n_i}$, where $1<n_1|n_2|\cdots|n_r$.
Then $r=r(\sG)$ is the {\it rank} of $\sG$.
It is natural to guess that $D(\sG)=\sum_{i=1}^rn_i-r+1=1+\sum_{i=1}^r(n_i-1)$
from pigeonhole intuition.
This was conjectured in \cite{O} and is true by Theorem \ref{cyclic} 
for rank 1 groups.
Moreover, it was proven in \cite{O} for rank 2 groups and $p$-groups as well.
However, it was proven in \cite{EBK,GS} that the conjecture is false 
for some groups each rank at least 4.
What remains open is the instance of rank 3.

\begin{conjecture}\label{davenport}
If $\sG$ is a finite abelian group of rank $r(\sG)=3$ then its
Davenport constant satisfies $D(\sG)=\sum_{i=1}^rn_i-r+1$.
\end{conjecture}

% ##########################################################################
%
% 		Pebbling Numbers
%

\section{Pebbling Numbers}\label{PebNum}

There are many known results regarding $\pi(G)$.
If one pebble is placed at each vertex other than the root vertex, $r$,
then no pebble can be moved to $r$.
Also, if $w$ is at distance $d$ from $r$, and $2^{d}-1$ pebbles are placed
at $w$, then no pebble can be moved to $r$.
Thus we have that $\pi(G) \ge max\{n(G), 2^{\diam(G)}\}$.
Graphs $G$ that satisfy $\pi(G)=n(G)$ are known as {\it Class 0} graphs, 
which include the complete graph $K_n$, the $d$-dimensional cube $Q^d$
\cite{Chung}, complete bipartite graphs $K_{m,m}$ \cite{Clarke}, 
and many others.
We will say more about such graphs in Section \ref{DCC}.
Any graph $G$ with a cut vertex $x$ has $\pi(G)>n(G)$.
(Indeed, let $v\in G_1$ and $u\in G_2$, where $G_1$ and $G_2$ are two
components of $G-x$.
Define the configuration $C$ by $C(v)=C(x)=0$, $C(u)=3$ and $C(w)=1$
for every other vertex $w$.
Then $|C|=n$ and $C$ cannot reach $v$.)
The path $P_n$, the cube $Q^d$ \cite{Chung}, the Petersen graph $P$
\cite{Clarke}, the even cycle $C_{2d}$ \cite{PSV}, and the line graph $L_n$
of the complete graph $K_n$ \cite{WW} are examples of graphs $G$ that 
satisfy $\pi(G)=2^{\diam(G)}$, while the odd cycle $C_{2d+1}$ \cite{PSV} 
is an example of a graph not satisfying either lower bound.
Another standard result is the pebbling number of a tree, which is worked 
out in \cite{Moews}.

Regarding upper bounds, it follows immediately from the Pigeonhole principle 
that a graph $G$ on $n$ vertices with diameter $d$ has pebbling number 
$\pi(G)\le (n-1)(2^d-1)+1$.
It would be interesting to find better general bounds on $\pi(G)$, especially
not involving $n$.
The independence number seems not to be useful.
For example, there is no function $g$ such that every graph $G$ of independence
number $\a$ and diameter $d$ has pebbling number $\pi(G)\le g(\a )2^d$.
Indeed, we define a family of graphs $G_m$ which satisfy $\diam(G)=d$ and
$\a (G)=2^{d-2}+1$, but which have pebbling number $\pi(G_m)\rar\infty$
as $m\rar\infty$.
Define $G_m=Q^{d-1}\cup K_m\cup E$, where $x\in V(Q^n)$ and the edge set
$E=\{xv\ |\ v\in V(K_m)\}$.
Since $x$ is a cut vertex we know by the above comment that
$\pi(G_m)>2^{d-1}+m$.
However, the domination number (minimum size of a dominating set)
can be useful.
Chan and Godbole \cite{CG} made the following improvements on the general
upper bound.

\begin{theorem}\label{betterupper}
Let $\dom(G)$ denote the domination number of $G$.
Then
\begin{enumerate}
\item
$\pi(G)\le (n-d)(2^d-1)+1$,
\item
$\pi(G)\le (n+\lfloor\frac{n-1}{d}\rfloor-1)2^{d-1}-n+2$, and
\item
$\pi(G)\le 2^{d-1}(n+2\dom(G))-\dom(G)+1$.
\end{enumerate}
The inequalities in parts 1 and 2 are sharp, and the coefficient of $2$ 
in part 3 can be reduced to $1$ in the case of perfect domination.
\end{theorem}

For any two graphs $G_1$ and $G_2$, we define the {\it cartesian product}
$G_1 \gbox G_2$ to be the graph with vertex set
$V(G_1\gbox G_2)=\{(v_1,v_2) | v_1 \in V(G_1), v_2 \in V(G_2)\}$ and edge set
$E(G_1\gbox G_2)=\{\{(v_1,v_2),(w_1,w_2)\} | (v_1=w_1$ and 
$(v_2,w_2)\in E(G_2))$
or $(v_2=w_2$ and $(v_1,w_1)\in E(G_1))$\}.
The cube $Q^d$ can be built recursively from the cartesian product, and
Chung's result \cite{Chung} that $\pi(Q^d)=2^d$ (and more generally 
Theorem \ref{PathProd} below) would follow easily from Graham's conjecture, 
which has generated a great deal of interest.

\begin{conjecture}(Graham)\label{Graham}
$\pi(G_1\gbox G_2)\le \pi(G_1)\pi(G_2)$.
\end{conjecture}

It is worth mentioning that there are some results which verify Graham's 
conjecture.
Among these, the conjecture holds for a tree by a tree \cite{Moews},
a cycle by a cycle \cite{Her,HH,PSV}, and a clique by a graph with the 
$2$-pebbling property
\cite{Chung} (see below).
Recently, Feng and Kim verified the conjecture for complete bipartite
graphs \cite{FK1} and for wheels or fans \cite{FK2}.
It is also proven in \cite{Chung} that the conjecture holds when
each $G_i$ is a path.
Let $P_n$ be a path with $n$ vertices and for ${\bf d}=\la d_1,\ldots,d_m\ra$
let $P_{\bf d}$ denote the graph $P_{d_1+1}\gbox\cdots\gbox P_{d_m+1}$.

\begin{theorem}\label{PathProd}
For nonnegative integers $d_1,\ldots,d_m$,
$\pi(P_{\bf d})=2^{d_1+\ldots +d_m}$.
\end{theorem}

The conjecture was also verified recently \cite{CH3} for graphs of high 
minimum degree, using Theorem \ref{Conn}.

\begin{result}\label{MinDeg}
If $G_1$ and $G_2$ are connected graphs on $n$ vertices that satisfy 
$\d(G_i)\ge k$ and $k\ge 2^{12n/k+15}$, then
$\pi(G_1\gbox G_2)\le \pi(G_1)\pi(G_2)$.
\end{result}
In particular, there is a constant $c$ so that if $k>cn/\lg{n}$ then
$G_1\gbox G_2$ is Class 0.
We will present probabilistic versions of Conjecture \ref{Graham} in Sections
\ref{GraphThresh} and \ref{PebblingThresh}.

\begin{figure}
\centerline{\includegraphics[height=1.5in]{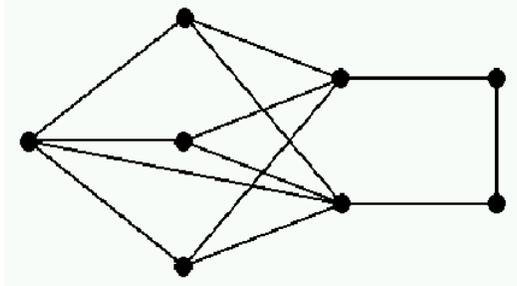}}
\caption{The Lemke graph}\label{lemke}
\end{figure}

A graph $G$ is said to have the {\it 2-pebbling property} if two pebbles can 
be moved to any specified vertex when the initial configuration $C$ has size
$2\pi(G)-s(C)+1$, where $s(C)$ is the number of vertices $v$ with $C(v)\ge 1$.
This property is crucial in the proof of Theorem \ref{PathProd}.
Until recently, only one graph $L$ (see Figure \ref{lemke}),
due to Lemke, was known not to have the
2-pebbling property, although a family of related graphs were conjectured in
\cite{Foster} not to have the property either.
Although that conjecture is still unresolved, Wang \cite{Wa} proved that 
a slight modification of Snevily's graphs yield an infinite family of 
graphs, none of which have the 2-pebbling property.
Also found in \cite{PSV} is the conjecture that all bipartite graphs
have the 2-pebbling property.
It is possible that the square of the Lemke graph might be a 
counterexample to Graham's Conjecture.

\begin{problem}
Find $\pi(L^2)$.
\end{problem}

In order to describe the generalization employed by Chung to prove the
Kleitman-Lemke theorem, we need to introduce generalized pebbling.
A {\it $q$-pebbling step} in $G$ consists of removing $q$
pebbles from a vertex $u$, and placing one pebble on a neighbor $v$ of $u$.
We say that a pebbling step from $u$ to $v$ is {\it greedy} ({\it semigreedy})
if $\dist(v,r)<\dist(u,r)$ ($\le$), where $r$ is the root vertex, and that a 
graph $G$ is {\it (semi) greedy} 
if for any configuration of $\pi(G)$ pebbles on the vertices of $G$ we can move 
a pebble to any specified root vertex $r$, in such a way that each pebbling 
step is (semi) greedy.

Let $P_{\bf d}=P_{d_1+1}\gbox\cdots\gbox P_{d_m+1}$ be a product of paths, 
where ${\bf d}=(d_1,\ldots,d_m)$.
Then each vertex $v\in V(P_{\bf d})$ can be represented by a vector
${\bf v}=\la v_1,\ldots,v_m\ra$, with $0\le v_i\le d_i$ for each $i\le m$.
Let ${\bf e}_i=\la 0,\ldots,1,\ldots,0\ra$, be the $i^{\it th}$ 
standard basis vector.
Denote the vector $\la 0,\ldots,0\ra$ by {\bf 0}.
Then two vertices $u$ and $v$ are adjacent in $P_{\bf d}$ if and only if 
$|{\bf u}-{\bf v}|={\bf e}_i$ for some integer $1\le i\le m$.
If ${\bf q}=(q_1,\ldots,q_m)$, then we may define {\bf q}-{\it pebbling} 
in $P_{\bf d}$ to be such that each pebbling step from {\bf u} to {\bf v} 
is a $q_i$-pebbling step whenever $|{\bf u}-{\bf v}|={\bf e}_i$.
We denote the {\bf q}-pebbling number of $P_{\bf d}$ by $\pi_{\bf q}(P_{\bf d})$.
Chungs's proof of Theorem \ref{PathProd} uses the following theorem \cite{Chung}. 
For integers $q_i,d_i\ge 1$, $1\le i\le m$, we use ${\bf q}^{\bf d}$
as shorthand for the product $q^{d_1}_1\cdots q^{d_m}_m$.

\begin{theorem}\label{GenPathProd}
Suppose that ${\bf q}^{\bf d}$ pebbles are assigned to the vertices of
$P_{\bf d}$ and that the root $r={\bf 0}$.
Then it is possible to move one pebble to $r$ via greedy {\bf q}-pebbling.
\end{theorem}
In addition, it was shown in \cite{CHH} that
$\pi_{\bf q}(P_{\bf d})={\bf q}^{\bf d}$, and moreover, that $P_{\bf d}$ is greedy.
Also in \cite{CHH} the following generalization of Graham's conjecture to
{\bf q}-pebbling was made.

\begin{conjecture}\label{GenGraham}
$f_{\bf q}(G_1\gbox G_2)\le f_{q_1}(G_1)f_{q_2}(G_2)$.
\end{conjecture}

Can Conjecture \ref{Graham} be proved in the case when both 
graphs share certain extra properties such as greediness and tree-solvability?
We say that a graph $G$ is {\it tree-solvable} if, whenever $C$ is a
configuration on $G$ of size $\pi(G)$, it is possible to solve $C$ in such
a way that the edges traversed by pebbling steps form an acyclic graph.
There are graphs which are neither greedy nor tree-solvable.
For example, let $\{a,b,c,d,e,f\}$ be the vertices (in order) of the cycle 
$C_6$, and form $G$ from $C_6$ by adjoining new vertices $g$ to $a$ and $c$, 
and $h$ to $a$ and $e$.
It is not difficult to show that $\pi(G)=9$.
With root $r=d$ and $C(a,b,f,g,h)=(1,3,3,1,1)$, $C$ has no greedy solution
(nor does it have a semi-greedy solution).
With $r=f$ and $C(b,c,d,g,h)=(1,5,1,1,1)$, $C$ has no tree-solution.
(It is worth noting that $G$ is bipartite.)

Another illustrative graph is $G=P_3\gbox S_4$, which is also bipartite.
($S_n$ is the {\it star} with $n$ vertices, also denoted $K_{1,n-1}$).
Easily, $\pi(P_3)=4$ and $F(S_4)=5$ (these numbers also follow from our
classification of diameter two graphs below).
Although each pebbling number is only one more than the number of vertices,
$\pi(G)=18$ is far greater than $n(G)=12$; it would be interesting to 
discover how much greater this gap can be for other graphs.
Also, notice that $18<\pi(P_3)\pi(S_4)$, a strict inequality.
More importantly, as observed by Moews \cite{M}, $G$ is not semi-greedy.
Indeed, think of $G$ as three pages of a book, let $r$ be the corner vertex
of one of the pages, $x$ the farthest corner vertex of another page, 
$u$, $v$ and $w$ the three vertices of the third page, and let
$C(u,v,w,x)=(1,1,1,15)$.
This shows that even semi-greediness is not preserved by the product of two
greedy graphs.
Probably, the following holds.

\begin{conjecture}\label{greedyconj}
Almost every graph is greedy and tree-solvable.
\end{conjecture} 

% ##########################################################################
%
% 		Diameter, Connectivity, and Class 0
%

\section{Diameter, Connectivity and Class 0}\label{DCC}

Because of the lower bound of $n(G)$ on the pebbling number $\pi(G)$,
it is natural to try to classify Class 0 graphs, or at least give conditions 
which either guarantee or prohibit Class 0.
This is of course extremely difficult, although some preliminary results
have proved quite interesting.
As argued above, graphs of connectivity 1 are not Class 0.
In \cite{PSV} we find the following theorem.

\begin{theorem}\label{Diam2}
If $\diam(G)=2$ then $\pi(G)=n(G)$ or $n(G)+1$.
\end{theorem}

Class 0 graphs of diameter 2 are classified in \cite{CHH}.
A particularly crucial graph in the characterization is the graph $G$, built
from the bipartite graph $C_6$ by connecting all the vertices 
of one of the parts of the bipartition to each other.
$G$ has connectivity 2 and diameter 2 but $\pi(G)>6=n(G)$, as witnessed by 
the configuration $C(a,b)=(3,3)$, where $a$, $b$ and the root $r$ are 
independent.
The following corollary to the characterization appears in \cite{CHH}.
Denote the connectivity of a graph $G$ by $\k(G)$.

\begin{result}\label{3Conn}
If $\diam(G)=2$, and $\k (G)\ge 3$ then $G$ is of Class 0.
\end{result}

From this it follows that almost all graphs (in the probabilistic sense)
are of Class 0, since almost every graph is 3-connected with diameter 2.
The following result, conjectured in \cite{CHH}, was proved in \cite{CHKT}.
This result was used to prove a number of other theorems, including
Theorems \ref{MinDeg}, \ref{Girth} and \ref{Class}.

\begin{result}\label{Conn}
There is a function $k(d)$ such that if $G$ is a graph with $\diam(G)=d$ and
$\k (G)\ge k(d)$ then $G$ is of Class 0.
Moreover, $2^d/d\le k(d)\le 2^{2d+3}$.
\end{result}

An upper bound for diameter 3 graphs was recently obtained in \cite{Buk}.

\begin{theorem}\label{Diam3}
If $\diam(G)=3$ then $\pi(G)\le 3n/2$, which is best possible.
\end{theorem}

Another guarantee for Class 0 membership may reside in the following, 
as yet unexplored, question.

\begin{question}\label{Greedy}
Is it true that every greedy graph is of Class 0?
\end{question}

A nice family of graphs in relation to Theorem \ref{Conn} is the following.
For $m\ge 2t+1$, the {\it Kneser graph}, $K(m,t)$, is the graph with
vertices ${[m] \choose t}$ and edges $\{ A,B\}$ whenever $A\cap B=\mt$.
The case $t=1$ yields the complete graph $K_m$ and the case $m=5$ and $t=2$
yields the Petersen graph $P$, both of which are Class 0.
When $t\ge 2$ and $m\ge 3t-1$ we have $\diam(K(m,t))=2$.
Also, it is not difficult to show that $\k(K(m,t))\ge 3$ in this range,
implying that $K(m,t)$ is Class 0 by Theorem \ref{3Conn}.
Furthermore, Chen and Lih \cite{Chen} have shown that $K(m,t)$ is
connected, edge transitive, and regular of degree ${m-t\choose t}$.
A theorem of Lov\'asz \cite{Lov} states that such a graph has connectivity
equal to its degree, and thus $\k=\k(K(m,t))={m-t\choose t}$.
Therefore, using Theorem \ref{Conn}, it is not difficult to prove
the following.

\begin{result}\cite{H}\label{Kneser}
For any constant $c>0$, there is an integer $t_0$ such that, for $t>t_0$,
$s\ge c(t/\log_2 t)^{1/2}$ and $m=2t+s$, we have that $K(m,t)$ is Class 0.
\end{result}

In the context of graph pebbling, the family of Kneser graphs is interesting
precisely because the graphs become more sparse as $m$ decreases toward
$2t+1$, so the diameter (as well as the girth) increases and yet the
connectivity decreases.

\begin{question}\label{ClassKneser}
For $1\le s\ll (t/\log_2 t)^{1/2}$, is $K(2t+s,t)$ Class 0?
\end{question}
Because Pachter et al. \cite{PSV} also proved that diameter two
graphs have the 2-pebbling property, it is interesting as well
to ask whether $K(m,t)$ has the 2-pebbling property when $m<3t-1$ (i.e.,
when its diameter is at least 3).
In addition, one should ask whether graphs of large enough connectivity
(in terms of the diameter) have the 2-pebbling property.
This would follow from Theorem \ref{Conn} if all Class 0 graphs have
the 2-pebbling property.

Regarding conditions which prohibit Class 0 membership, one can easily show
that if $\gir(G)>2\log n$ then $\pi(G)>n(G)$.
The following question was asked in \cite{H}:
Is there a constant $g$ so that $\gir(G)>g$ implies $\pi(G)>n(G)$?
This question was answered in the negative in \cite{CH3} by using Theorem
\ref{Conn}, along with a probabilistic (deletion) method analogous to Erd\H os's
construction \cite{E} of graphs of arbitrarily high girth and chromatic number.

\begin{result}\label{Girth}
Let $g_0(n)$ denote the maximum number $g$ such that there exists a 
Class 0 graph $G$ on at most $n$ vertices with finite $\gir(G)\ge g$.
The for all $n\ge 3$ we have
$$\lfloor\sqrt{(\lg n)/2+1/4}-1/2\rfloor\le g_0(n)\le 1+2\lg n\ .$$
\end{result}

Finally, denote by $B(m,k)$ the family of all connected, bipartite graphs on 
$n=2m$ vertices with minimum degree at least $k$.
What can be said of the pebbling numbers of graphs in $B(m,k)$?
In particular, the only graph in $B(m,2)$ is $C_{2m}$, having pebbling
number $2^m$.
The only graph in $B(m,m)$ is $K_{m,m}$, having pebbling number $n=2m$,
as first shown in \cite{Clarke}.
Also shown in \cite{Clarke} is the result that if $G\in B(m,m-1)$ then
$G$ is of Class 0.
In \cite{Hurl} we find that if $k\ge 2m/3$ (and $m\ge 6$) then every
$G\in B(m,k)$ is of Class 0.
Define $b(m)$ to be the smallest $k$ such that $G$ is of Class 0 for 
every $G\in B(m,k)$.
Recently in \cite{CH1} was proved the following result.

\begin{result}\label{RegularBipartite}
For all $m\ge 7$ we have $b(m)\ge\lfloor m/2\rfloor +1$.
For all $m\ge 336$ we have $b(m)=\lfloor m/2\rfloor +1$.
\end{result}
If $G$ has some pseudo-random properties then the bound for $t$ seems
to be significantly smaller. 
As far as sparse graphs are concerned, a related problem is to find
the transitions for $k$ at which the pebbling numbers of the graphs in
$B(m,k)$ change from exponential to subexponential in $n$, and from
superlinear to linear in $n$.  

In connection with Theorem \ref{RegularBipartite}, it is natural to
consider the family $G(n,k)$ of connected, minimum degree $k$ graphs instead of
bipartite graphs. 
In \cite{CH1} is proved an analogous result.
Let $g(n)$ be the minimum $k$ such that every graph in $G(n,k)$ is Class 0.

\begin{result}\label{Normal} 
For all $n\ge 6$ we have $g(n)\le\lfloor n/2\rfloor$.
For all $n\ge 9$ we have $g(n)=\lfloor n/2\rfloor$.
\end{result}

% ##########################################################################
%
% 		Cover Pebbling and Optimal Pebbling
%

\section{Cover Pebbling and Optimal Pebbling}\label{CoverPebbOpt}

In terms of possible applications it may be relevant to consider having to
move a pebble to every vertex simultaneously, rather than to just one root.
More generally, given a weight function $w$ on the vertices of a graph $G$,
We define the {\it weighted cover pebbling number} $\c_w(G)$ to be the
minimum number of pebbles needed to place, after a sequence of pebbling steps, 
$w(v)$ pebbles on vertex $v$, for all $v$, regardless of the initial 
configuration.
The {\it cover pebbling number} $\c(G)$ regards the case $w(v)=1$ for all $v$.
We say that $w$ is {\it positive} if $w(v)>0$ for every vertex $v$.
In \cite{CCFHPST}, for positive weight functions the weighted cover pebbling 
number for cliques and trees is found, and it is shown that the ratio between 
this parameter and the pebbling number can be arbtrarily large, even within 
the class of trees.
Given $w$ define $|w|=\sum_v w(v)$ and $\min w=\min_v w(v)$.
The results are as follows.

\begin{result}\label{wCovPebbK}
For every positive weight function $w$ we have $\c_w(K_n)=2|w|-\min w$.
\end{result}

The {\it fuse} $F_l(n)$ is the graph on $n=m+k$ vertices composed of a path
(or {\it wick}) on $l$ vertices with $n-l$ independent vertices ({\it sparks})
incident to one of its endpoints.

\begin{result}\label{wCovRat}
For every $n$ and $l$ we have $\c(F_l(n))=(n-l+1)2^l-1$.
Thus, for $n=2^l+l$, we have $\c(F_l(n))/\pi(F_l(n))>(n-\lg n)/2$.
\end{result}

Given a graph $G$ define $s_w(G)=\max_v\sum_u w(u)2^{\dist(u,v)}$, where 
$\dist(u,v)$ denotes the distance between $u$ and $v$.

\begin{result}\label{wCovPebbT}
For every positive weight function $w$ we have $\c_w(T)=s_w(T)$
for every tree $T$.
\end{result}

It would be interesting to find the weighted pebbling numbers of other graphs 
such as cycles, cubes, complete $r$-partite graphs and random graphs, among 
others.
The main technique in proving Theorem \ref{wCovPebbT} involves showing that 
the largest $w$-unsolvable configuration is {\it simple}; that is, 
concentrated on a single vertex.
The same cannot be said in general for weight functions that are not positive.
An important question to investigate along these lines is whether, for 
arbitrary graphs and positive weight functions, the largest $w$-unsolvable 
configuration is always simple.
A positive answer verifies the so-called Stacking Conjecture (now Theorem), 
below.

\begin{theorem}\label{stack}
For every graph $G$ and positive weight function $w$ we have $\c_w(G)=s_w(G)$.
\end{theorem}

This theorem originally was verified in \cite{HM1} for the case of covering 
$d$-dimensional cubes ($\c(Q^n)=3^n$).
Other graphs were considered in \cite{TW,WY}.
Theorem \ref{stack} was recently proved in general independently in
\cite{S,VW}.

One should also search for classes of graphs for which there is some constant 
$c$ so that every graph $G$ in the class satisfies $\c(G)/\pi(G)\le c$.
For example, the classes of cliques and paths both have $c=2$.

Finally, we define the weighted pebbling number $\pi_\bw(G)$ to be the minimum 
number $t$ so that, for every weight function $w$ with $|w|=\bw$, 
every configuration of $t$ pebbles is $w$-solvable.
If one is concerned only with positive $w$, then such an evaluation on any 
tree $T$ is trivial.
One would compute 
$\max_w\max_v\sum_u2^{d(u,v)}
= \max_v\max_w\sum_u2^{d(u,v)}
= \max_v |w|2^{\max_ud(u,v)}
= \w2^{\diam(T)}$.
However, nonpositive $w$ must be considered, and thus finding weighted 
cover pebbling numbers for such $w$ is critical.

Instead of looking for the size of the largest unsolvable configuration 
in a graph, which is essentially the task in finding the pebbling number 
of the graph (minus one), one could look for the size of the smallest 
solvable configuration, which is the task in finding what is called the 
{\it optimal} pebbling number, $\pi^*(G)$.
The first results in this direction showed that $\pi^*(P_n)=\lceil 2n/3\rceil$
\cite{PSV}.
In \cite{BCCMW} it is proved that $\pi^*(G)\le\lceil 2n/3\rceil$ for all $G$
(equality also holds for cycles),
and the lower bound $\pi^*(Q^d)\ge (4/3)^d$ is found in \cite{M2}.
Caterpillars, cycles and other graphs have been considered in \cite{FW,FS1,FS2},
and the following interesting analog of Graham's conjecture was recently
proven in \cite{FW,FS1}.

\begin{theorem}\label{OptGraham}
For all graphs $G$ and $H$ we have $\pi^*(G\gbox H)\le \pi^*(G)\pi^*(H)$.
\end{theorem}

What is most surprising is that, while $\pi(Q^n)$ is known exactly, only
$\pi^*(Q^n)=(\frac{4}{3})^{n+O(\log n)}$ is known at present \cite{M2}.
Clearly, much work needs to be done in this area.
Some of the newest results consider graphs of high minimum degree.
The following result of Czygrinow appears in \cite{BCCMW}.

\begin{theorem}\label{OptDeg}
If $G$ is a connected graph with $n$ vertices and $\d(G)=k$, then
$\pi^*(G)\le\frac{4n}{k+1}$.
\end{theorem}
The result is not known to be sharp, as given by the following result in
\cite{BCCMW}.

\begin{theorem}\label{OptDegLow}
For all $t\ge 1$, $k=3t$ and $n\ge k+3$, there is a graph $G$ with $\d(G)=k$ and
$\pi^*(G)\ge (2.4-\frac{24}{5k+15}-o(1))\frac{n}{k+1}$.
\end{theorem}
The graphs discovered for this theorem are clever modifications of a blow-up 
of the vertices of a cycle into cliques.
It is noted that the $o(1)$ term disappears when $n$ is a multiple of $t+1$.
It would be interesting to discover the correct asymptotic coefficient of
$\frac{n}{k+1}$, somewhere between 2.4 and 4.
The authors also ask if the general upper bound of $\lceil 2n/3\rceil$ can
be improved because of high minimum degree.

\begin{theorem}\label{nOver2}
Is it true that $\pi^*(G)\le\lceil n/2\rceil$ for all graphs $G$ with
$\d(G)\ge 3$?
\end{theorem}
Such a bound would be sharp for $n(G)\ge 6$.

If girth is also considered then one can say more.
Let $c_k(t)=1+k\sum_{i=0}^{t-1}(k-1)^i$ and
$c^\pr(t)=(2^{2t}-2^{t+1})\frac{t}{t-1}$.
The following theorem of \cite{BCCMW} displays an asymptotic bound of $3n/8$.

\begin{theorem}\label{DegGir}
Let $k\ge 3$, $t\ge 2$ and $(k,t)\not\in(3,2)$.
Then every $n$-vertex graph $G$ with $\d(G)=k$ and $\gir(G)\ge 2t+1$
satisfies $\pi^*(G)\le 2^{2t}n/(c_k(t)+c^\pr(t))$.
\end{theorem}

The optimal pebbling numbers of linear ($P_m\gbox K_2$), circular ($C_m\gbox K_2$)
and M\" obius (circular with a twist) are also determined in \cite{BCCMW}:
$\pi^*=m$ unless $m\in\{2,5\}$, with $m$ as a lower bound always.

In terms of a taxi service, army unit, or other operational outfit being
able to strategically place their limited resources so as to reach any
desired location, this line of research may be most fruitful and significant. 

% ##########################################################################
%
% 		Fractional Pebbling and Complexity
%

\section{Fractional Pebbling and Complexity}\label{FracPebbComp}

Many combinatorial optimization parameters have fractional counterparts.
Graph parameters like matching number, chromatic number \cite{SU} are 
well known examples, as is the dimension of posets \cite{BS}, and the 
pebbling number is no exception.
Let $\pi(G,k)$ be the minimum $t$ such that, for every choice of root $r$ 
and for every configuration of $t$ pebbles on $G$, it is possible to move $k$ 
pebbles to $r$.
Then the fractional pebbling number $\hat{\pi}(G)=\lim_{k\to\infty}\pi(G,k)/k$.
For example, $\hat{\pi}(K_n)=\lim_{k\to\infty}(n+2k-2)/k=2$.
Also, using Moews's tree theorem \cite{Moews}, it is easy to compute
$\hat{\pi}(T)=2^{\diam(T)}$ for every tree $T$.
The following general result was recently proven in \cite{HK}.

\begin{result}\label{FracPebbThm}
For every graph $G$ we have that $\hat{\pi}(G)=2^{\diam(G)}$.
\end{result}

Typically, the parameter in question can be formulated in terms of an integer 
linear program, and the relaxation of the program gives rise to the fractional 
version of the parameter.
In the case of graph pebbling, while it is routine to verify in polynomial time
that a particular configuration reaches a particular root, no simple matrix 
condition exists to capture whether every configuration of a fixed size $t$ is 
solvable for a particular root.
Thus minimizing $t$ as an integer program remains elusive, and so it is 
difficult 
to say at present whether or not a parameter exists that is dual to $\pi$.
Recently, Lourdusamy and Somasundaram \cite{LS} found a proof that verifies 
Graham's conjecture for the case $C_5\gbox C_5$.
Their proof uses linear programming and lends credence to the above ideas 
while suggesting that further explorations into the connections between 
linear programming duality and pebbling numbers may be fruitful and 
worthwhile indeed.

Because it is easy to verify a pebbling solution, the question of deciding
whether a given configuration $C$ reaches a given root $r$ is in {\sf NP}.
In fact it has been shown \cite{HK} that this question is at least has hard as 
{\sf NP-complete} (see also \cite{MC,Wat}).

\begin{result}\label{Hard}
For every $4$-regular hypergraph $H$, there is a graph $G$, a configuration
$C$, and a root $r$, so that deciding if $C$ reaches $r$ in $G$ is at least
as hard as deciding if $H$ has a perfect matching.
\end{result}

A complexity upper bound is found in \cite{MC}, where it shown that
deciding whether $\pi(G)\le k$ is a $\Pi_2^P$-complete problem.
This means that it is complete for the class of languages computable
in polynomial time by coNP machines equipped with an oracle for an
NP-complete language.
They also proved that deciding whether $\pi^*(G)\le k$ is NP-complete,
while the same result is proved in \cite{Wat} for deciding whether
$\c(G)\le k$.

Instead of computing $\pi(G)$ exactly, one could try to approximate it. 
In other words the following question is of interest.

\begin{question}\label{comp2}
Does there exist a (constant/linear/polynomial) function $F(n)$ and an 
algorithm which, given a graph $G$ on $n$ vertices, can find a number in
polynomial time $\overline{\pi}(G)$ such that
$$|\pi(G) - \overline{\pi}(G)| \leq F(n)?$$
\end{question}

Just like in the case of other ``hard'' problems one can  design a ``good'' 
approximation algorithm for special cases of graphs. 
In particular, is it possible to design an approximation algorithm for the 
case when $G$ is a dense graph?
Perhaps an easier and more natural question is to ask for an algorithm which, 
given a configuration and a root, will pebble to the root efficiently.

\begin{question}\label{comp3}
Does there exist an algorithm which, given a graph $G$, a root $r$, and an 
$r$-solvable configuration $C$, can pebble to $r$ in time polynomial in $n$?
\end{question} 

Of course, Theorem \ref{Hard} shows that this is in general impossible 
unless ${\sf P}={\sf NP}$, although the question for appropriately restricted 
classes of graphs, such as cubes, remains interesting.
Note that for some graphs (like trees) the greedy algorithm
will pebble to $r$ efficiently. 
However, as noted in Section \ref{PebNum}, there are graphs which are
neither tree-solvable, nor greedy (nor semi-greedy), and so the
question of how to pebble is far from obvious.
Even in the case of greedy graphs like $Q^d$, there are solvable configurations
that require decision-making; that is, some greedy approaches fail while
others succeed.

% ##########################################################################
%
% 		Graph Thresholds
%

\section{Graph Thresholds}\label{GraphThresh}

The notion that graphs with very few edges tend to have large pebbling
number and graphs with very many edges tend to have small pebbling
number can be made precise as follows.
Let $\gnp$ be the random graph model in which each of the $n\choose 2$ 
possible edges of a random graph having $n$ vertices appears independently 
with probability $p$.
For functions $f$ and $g$ on the natural numbers we write that
$f\ll g$ (or $g\gg f$) when $f/g\rar 0$ as $n\rar\infty$.
Let $o(g)=\{f\ |\ f\ll g\}$ and define $O(g)$ (resp., $\Omega(g)$) to be 
the set of functions $f$ for which there are constants $c, N$ such that
$f(n)\le cg(n)$ (resp., $f(n)\ge cg(n)$) whenever $n>N$.
Finally, let $\Th(g)=O(g)\cap\Omega(g)$.

Let $\cgQ$ be a property of graphs and consider the probability
$\Pr(\cgQ)$ that the random graph $\gnp$ has $\cgQ$.
For large $p$ it may be that $\Pr(\cgQ)\rar 1$ as 
$n\rar\infty$, and for small $p$ it may be that 
$\Pr(\cgQ)\rar 0$ as $n\rar\infty$.
More precisely, define the {\it threshold} of $\cgQ$, $\T(\cgQ)$, 
to be the set of functions $t$ for which $p\gg t$ implies that
$\Pr(\cgQ)\rar 1$ as $n\rar\infty$, and $p\ll t$ implies that
$\Pr(\cgQ)\rar 0$ as $n\rar\infty$.

It is not clear that such thresholds exist for arbitrary $\cgQ$.
However, we observe that Class 0 is a monotone property (adding edges
to a Class 0 graph maintains the property), and a theorem of Bollob\'as
and Thomason \cite{BT} states that $\T(\cgQ)$ is nonempty for every 
monotone $\cgQ$. 
It is well known \cite{ER1} that $\T({\rm connected})=\Th(\lgn /n)$, 
and since connectedness is required for Class 0, we see that 
$\T({\rm Class\ 0})\subseteq\Omega(\lgn /n)$.
In \cite{CHH} it is noted that $\T({\rm Class\ }0)\subseteq O(1)$.
In \cite{CHKT} Theorem \ref{Conn} was used to prove the following result.

\begin{result}\label{Class}
For all $d>0$, $\T({\rm Class\ }0)\subseteq O((n\lgn)^{1/d}/n)$.
\end{result}

As a probabilistic version of Theorem \ref{RegularBipartite}, we might
replace the probability space of $\gnp$ by $\bmp$, the random graph model
in which each of the $m^2$ possible edges of a random bipartite graph
having $m$ vertices in each bipartition appear independently with
probability $p$.
One might be interested in the following problem.

\begin{problem}\label{RandBip}
Find the Class 0 threshold $\T_{\B}({\rm Class\ }0)$ of the random
bipartite graph $\bmp$.
\end{problem}

We now reconsider Graham's Conjecture \ref{Graham} from a probabilistic
point of view.
We would like to say at the very least that it holds almost always.
More precisely, we offer the following.

\begin{conjecture}\label{RandGraham}
Let $G,H\in\gnp$ for any $p=p(n)$.
Then $\Pr [\pi(G\gbox H)\le \pi(G)\pi(H)]\rar 1$ as $n\rar\infty$.
\end{conjecture}

If it turns out that the threshold for Class 0 is bigger than that for
connectivity, then it would be a very interesting to investigate
the pebbling numbers of graphs in $\gnp$ for $p$ in that range.

Another important problem, especially in light of the previous
conjecture, is the following.

\begin{problem}\label{2PPThresh}
Find the threshold $\T(2{\rm PP})$ for the 2-pebbling property of the
random graph $\gnp$.
\end{problem}

% ##########################################################################
%
% 		Pebbling Thresholds
%

\section{Pebbling Thresholds}\label{PebblingThresh}

For this section we will fix notation as follows.
The vertex set for any graph on $N$ vertices will be taken to be
$\{v_i\ |\ i\in [N]\}$, where $[N]=\{0,\ldots,N-1\}$.
That way, any configuration $C:V(G_n)\rar\N$ is independent of $G_n$.
(Here we make the distinction that $n$ is the index of the graph $G_n$ in a 
sequence $\G=(G_1,\ldots,G_n,\ldots)$, whereas $N=N(G_n)$ denotes its number 
of vertices.)
Let $\K=(K_1,\ldots,K_n,\ldots)$ denote the sequence of complete graphs,
$\P=(P_1,\ldots,P_n,\ldots)$ the sequence of paths, and
$\Q=(Q^1,\ldots,Q^n,\ldots)$ the sequence of $n$-dimensional cubes.
Let $C_n:[N]\rar \N$ denote a configuration on $N=N(G_n)$ vertices.

Let $h:\N\rar\N$ and for fixed $N$ consider the probability space $X_N$
of all configurations $C_n$ of size $h=h(N)$, 
We denote by $P_N^+$ the probability that $C_N$ is $G_n$-solvable and
let $t:\N\rar\N$.
We say that $t$ is a {\it pebbling threshold} for $\G$, and write
$\t(\G)=\Th(t)$, if $P_N^+\rar 0$ whenever $h(N)\ll t(N)$ and $P_N^+\rar 1$
whenever $h(N)\gg t(N)$.  
The existence of such thresholds was recently established in \cite{BBCH}.

\begin{result}\label{ExistThresh}
Every graph sequence $\G$ has nonempty $\t(\G)$.
\end{result}

The first threshold result is found in \cite{Clarke}.
The result is merely an unlabeled version of the so-called 
``Birthday problem'', in which one finds the probability that 2 of $t$ 
people share the same birthday, assuming $N$ days in a year.

\begin{result}\label{CliqueThresh}
$\t(\K)=\Th(\sqrt N)$.
\end{result}
The same threshold applies to the sequence of stars ($K_{1,n}$).

It was discovered in \cite{CEHK} that every graph sequence $\G$ satisfies
$\t(\K)\lesssim\t(\G)\lesssim\t(\P)$, where $A\lesssim B$ is meant to 
signify that $a\in O(b)$ for every $a\in A, b\in B$.
The authors also discovered that $\t(\G)\sub O(N)$ when $\G$ is a sequence 
of graphs of bounded diameter, that $\t(\Q)\sub O(N)$, and that 
$\t(\P)\sub \W(N)$.  
Surprisingly, the threshold for the sequence of paths has not been
determined.
The lower bound found in \cite{CEHK} was improved in \cite{BBCH} to 
$\t(\P)\sub\W(N2^{c\slN})$ for every $c<1/\sqrt 2$, while the upper bound 
of $\t(\P)\sub O(N2^{2\slN})$ found in \cite{BBCH} was improved in 
\cite{GJSW} to $\t(\P)\sub O(N2^{c\slN})$ for every $c>1$.
Finally the lower bound was tightened recently in \cite{CH4} to nearly
match the upper bound.

\begin{result}\label{lower}
For any constant $c<1$, we have $\t(\P)\sub\W(N2^{c\slN})$.
\end{result}

Frustratingly, this still leaves room for a wide range of possible
threshold functions.
It is interesting that even within the family of trees, the pebbling
thresholds can vary so dramatically, as in the case for paths and stars.
Diameter seems to be a critical parameter.
It is quite natural to guess that families of graphs with higher
pebbling numbers have a higher threshold, but this kind of monotonicity
result remains unproven.

\begin{conjecture}\label{MonoThresh}
If $\pi(G_n)\le \pi(H_n)$ for all $n$ then $\t(\G)\le \t(\cH)$.
\end{conjecture}
In fact, this conjecture remains unproven even in the case that $\G$ and
$\cH$ are sequences of trees.
Moreover, there is some reason to believe the conjecture may be false, 
since the pebbling number is a worst-case scenario, while the threshold
is an average-case scenario.

Consider the following.
For a positive integer $t$ and a graph $G$ denote by $p(G,t)$ the probability
that a randomly chosen configuration $C$ of size $t$ on $G$ solvable.
Then Conjecture \ref{MonoThresh}
would follow from the statement that, if $\pi(G_n)\le \pi(H_n)$ then for
all $t$ we have $p(G_n,t)\ge p(H_n,t)$.
Unfortunately, although seemingly intuitive, this implication is false.
Using the Class 0 pebbling characterization theorem of \cite{CHH}, 
in \cite{CEHK} is found a family of pairs of graphs $(G_n,H_n)$, one pair
for each $n=3k+4$, for which the implication fails.
However, the implication may yet hold when $G_n$ and $H_n$ are trees.

\begin{question}\label{trees}
Let $G_n$ and $H_n$ both be trees on $n$ vertices so that 
$\pi(G_n)\le \pi(H_n)$.
Is it true that, for all $t$ (or those $t$ ``near'' one of the thresholds),
we have $p(G_n,t)\ge p(H_n,t)$?
\end{question}

At the heart of such investigations into monotonicity is the following,
most natural conjecture.

\begin{conjecture}\label{Spec}
Let $t_1,t_2$ be any functions satisfying
$\t(\K)\lesssim t_1\ll t_2\lesssim \t(\P)$.
Then there is some graph sequence $\G$ such that 
$t_1\lesssim\t(\G)\lesssim t_2$.
\end{conjecture}
This conjecture was proven in \cite{CH4} in the case that $\t(\P)$ is 
replaced by $\Th(N)$.
In fact, the family of fuses (defined in Section \ref{CoverPebbOpt}) 
covers this whole range.
What behavior lives above $\Th(N)$ remains unknown.

It is interesting to consider a pebbling threshold version of
Graham's conjecture.
Given graph sequences $\F$ and $\G$, define the sequence 
$\cH=\F\gbox\G=\{F_1\gbox G_1,\ldots,F_n\gbox G_n,\ldots\}$.
Suppose that $f(R)\in\t(\F)$, $g(S)\in\t(\G)$, and $h(T)\in\t(\cH)$, 
where $R=N(F_n)$, $S=N(G_n)$, and $T=N(H_n)=RS$.

\begin{question}\label{ProdThresh}
Is it true that, for $\F$, $\G$, and $\cH$ as defined above, we have
$h(T)\in O(f(R)g(S))$?
\end{question}
In particular, one can define the sequence of graphs $\G^k$ in the
obvious way.
In \cite{CH3} one finds tight enough bounds on $\t(\P^k)$ to show that
the answer to this question is yes for $\F=\P^i$ and $\G=\P^j$.
Another important instance is $\cH=\K^2$.
Boyle \cite{Boy} proved that $\t(\K^2)\in O(N^{3/4})$.
This was improved in \cite{BH}, answering Question \ref{ProdThresh} 
affirmatively for squares of complete graphs.

\begin{result}\label{SquareCliques}
For $\K^2=\{K_1^2,\ldots,K_n^2,\ldots\}$ we have $\t(\K^2)=\Th(N^{1/2})$.
\end{result}
This result is interesting because, by squaring, the graphs become fairly
sparse, and yet their structure maintains the low pebbling threshold.
The proof of the result tied the behavior of pebbling in $K_n^2$ to
the existence of large components in various models of random complete
bipartite graphs.

Another interesting related sequence to consider is
$\P_l=\{P_l^1,\ldots,P_l^n,\ldots\}$.
When $l=2$ we have $\P_2=\Q$, and the best result to date is the following
theorem of \cite{CW} (obtained independently in \cite{A}).

\begin{theorem}\label{Cubes}
For the sequence of cubes we have
$\t(\Q)\in\W(N^{1-\e})\cap O(N/(\lg\lg N)^{1-\e})$
for all $\e>0$.
\end{theorem}

Let $l=l(n)$ and $d=d(n)$ and denote by $\P_l^d$ the graph sequence
$\{P_{l(n)}^{d(n)}\}_{n=1}^{\infty}$, where $P_l^d=(P_l)^d$.
Most likely, fixed $l$ yields similar behavior to Theorem \ref{Cubes}.

\begin{conjecture}\label{Grids}
For fixed $l$ we have $\t(\P_l^n)\in o(N)$.
\end{conjecture}
In contrast, the results of \cite{CH3} show that $\t(\P^d)\in\W(N)$ for
fixed $d$.
Thus it is reasonable to believe that there should be some relationship 
between the two functions
$l=l(n)$ and $d=d(n)$, both of which tend to infinity, for which the 
sequence $\P_l^d$ has threshold on the order of $N$.

\begin{problem}\label{GridBalance}
Find a function $d=d(n)$ for which $\t(\P^d)=\Th(N)$.
In particular, how does $d$ compare to $n$?
\end{problem}

Finally one might consider the behavior of graphs of high minimum degree.
Define $\bG(n,\d)$ to be the set of all connected graphs on $n$ vertices
having minimum degree at least $\d=\d(n)$.
Let $\G_\d=\{G_1,\ldots,G_n,\ldots\}$ denote any sequence of graphs with
each $G_n\in\bG(n,\d)$.
In \cite{CH1} is proven the following.

\begin{result}\label{Dense} 
For every function $n^{1/2}\ll \d=\d(n)\le n-1$,
$\t(\G_\d)\sub O(n^{3/2}/\d)$.
In particular, if in addition $\d\in\W(n)$ then $\t(\G_\d)=\Th(n^{1/2})$.
\end{result}

% ##########################################################################
%
% 		Other Applications
%

\section{Other Variations and Applications}

A very recent incarnation of graph pebbling is introduced in \cite{GLP},
in which the {\it critical pebbling number} of a graph $G$ is defined.
A configuration $C$ of pebbles on $G$ is {\it minimally solvable} if it
is solvable but the removal of any pebble leaves it unsolvable.
While the optimal pebbling number measures the size of the smallest
minimally solvable configuration on $G$, the critical pebbling number
measures the size of the largest minimally solvable configuration.

Another variation is developed in \cite{GGTVWY}, combining ideas from
Sections \ref{PebNum} and \ref{CoverPebbOpt}.
The {\it support} of a configuration $C$ is the set of vertices $v$ having
$C(v)>0$.
The {\it domination cover pebbling number} $\psi(G)$ is the minimum number
$t$ so that from every configuration of $t$ pebbles one can reach a
configuration whose support is a dominating set of $G$.
Their motivation stems from transporting devices from initial positions
to eventual positions that allow them to monitor the entire graph.

Finally, graph pebbling investigations on directed graphs \cite{Gunda}
have begun as well.

Because of the superficial similarity of graph pebbling to other positional 
games on graphs, like ``Cops-and-Robbers'' \cite{NW,ST} and ``Chip-Firing'' 
\cite{ALSSTW,CE} for instance, and the high degree
of applicability of many of these games to structural graph theory \cite{RS}
and theoretical computer science \cite{ACP}, one shouldn't neglect the 
possibility that graph pebbling will have similar impact.
For example, one can think of the loss of a pebble during a pebbling step
as a toll or as a loss of information, fuel or electrical charge.
$q$-pebbling is one generalization of this rate of loss; another is simply
to choose any fixed rate $\a$ of loss.
In any case, instead of restricting the initial configuration to integer
values, let $C$ range among all nonnegative reals.
A pebbling step removes weight $x$ from one vertex and places weight 
$\a x$ at an adjacent vertex, for some fixed $0<\a <1$.
Still the objective is to place weight 1 at any prescribed root $r$ so that
there is enough money, fuel, information, or energy at that location in the 
network.
Of course, all of the questions raised herein may be asked about this
more general $\a$-{\it pebbling}.
It is conceivable that chip-firing may even come into play as a useful
model.
For a given graph $G$, one might be able to build an auxiliary graph $H$, 
so that chip-firing results on $H$ can be brought to bear on $\pi(G)$.
This opens up the theory to questions of the eigenvalues of the Laplacian
of $G$, and so on.

% ##########################################################################
% ##########################################################################
% 
%		BIBLIOGRAPHY
%

\end{document}